\documentclass[12pt]{article} %%% use \documentstyle for old LaTeX compilers

\usepackage[english]{babel} %%% 'french', 'german', 'spanish', 'danish', etc.
\usepackage{amssymb}
\usepackage{amsmath}
\usepackage{txfonts}
\usepackage{mathdots}
\usepackage[classicReIm]{kpfonts}
\usepackage[dvips]{graphicx} %%% use 'pdftex' instead of 'dvips' for PDF output
\usepackage{xcolor}
\usepackage[colorlinks]{hyperref}
\hypersetup{colorlinks=true,linkcolor=blue, citecolor=blue}

\usepackage[paper=a4paper, left = 2cm, right = 2cm, top = 2cm, bottom = 2cm]{geometry}
\begin{document}

%\selectlanguage{english} %%% remove comment delimiter ('%') and select language if required

\noindent 

\noindent 

\noindent 

\noindent 

\noindent 
\begin{center}
\section*{On higher order Codazzi tensors\\
on complete Riemannian manifolds}
\end{center}

\begin{center}
\noindent \textbf{    I. G. Shandra}${}^{1      }$\textbf{S. E. Stepanov}${}^{1,2}$

\noindent 

\begin{footnotesize}
\noindent ${}^{1 }$Dept. of Data Analysis and Financial Technology, Finance University, 

\noindent 49-55, Leningradsky Prospect, 125468 Moscow, Russia 

\noindent e-mail: ma-tematika@yandex.ru

\noindent 

\noindent ${}^{2}$ Dept. of Mathematics, Russian Institute for Scientific and 

\noindent Technical Information of the Russian Academy of Sciences, 

\noindent 20, Usievicha street, 125190 Moscow, Russia

\noindent e-mail: s.e.stepanov@mail.ru

\end{footnotesize}
\end{center}
\noindent 

\noindent \textbf{}

\noindent \textbf{Abstract. }{\small We prove several Liouville-type non-existence theorems for higher order Codazzi tensors and classical Codazzi tensors on complete and compact Riemannian manifolds, in particular. These results will be obtained by using theorems of the connections between the geometry of a complete smooth manifold and the global behavior of its subharmonic functions. In conclusion, we show applications of this method for global geometry of a complete locally conformally flat Riemannian manifold with constant scalar curvature because its Ricci tensor is a Codazzi tensor and for global geometry of a complete hypersurface in a standard sphere because its second fundamental form is also a Codazzi tensor.}

\noindent 

\noindent \textbf{Keywords}:\textbf{ }complete Riemannian manifold, higher order Codazzi tensor, subharmonic function. 

\noindent \textbf{Mathematical Subject Classification: }53C20; 53C25; 53C40

\noindent              
\pagebreak

\noindent \textbf{1. Introduction}

\noindent Let $\left({\rm M,g}\right)$ be a Riemannian manifold of dimension \textit{n} $\geq$ 2 with the Levi-Civita connection $\nabla $. Everywhere in the following we denote by $\Lambda ^{{\rm q}} M$ and $S^{p} M$ the vector bundles of exterior differential $q$-forms $\left(1\le q\le n-1\right)$ and symmetric differential $p$-forms $\left(p\ge 2\right)$ on ${\rm M}$. Throughout this paper we will consider the vector spaces of their \textit{C}${}^{\infty}$-sections denoted by \textit{${\rm C}^{\infty } \Lambda ^{{\rm p}} M$} and $C^{\infty } S^{p} M$, respectively. The Riemannian metric $g$ induces a point-wise inner product on the fibres of each of these spaces. 

\noindent A symmetric bilinear form ${\rm T}\in {\rm C}^{\infty } S^{2} M$on a Riemannian manifold $\left(M,\, g\right)$ is said to be a \textit{Codazzi tensor} if it satisfies the equation [1]; [2, p. 435]; [3, pp. 24; 56; 68]
\begin{equation} \label{GrindEQ__1_1_} 
\left(\nabla _{X} T\right)\, \left(Y,\, Z\right)=\left(\nabla _{Y} T\right)\, \left(X,\, Z\right) 
\end{equation} 
for any tangent vector fields $X,Y,Z$ on ${\rm M}$. The Codazzi tensor ${\rm T}\in {\rm C}^{\infty } S^{2} M$ is called \textit{non-trivial} if it is not a constant multiple of metric [1, p. 15]. Alongside with it we know from [1, p. 17] that every smooth manifold ${\rm M}$ carries a $C^{\infty } $-metric $g$ such that $\left(M,\, g\right)$ admits a nontrivial Codazzi tensor ${\rm T}\in {\rm C}^{\infty } S^{2} M$. We remark here that this result is essentially local.

\noindent Codazzi tensors appear in a natural way in many geometric situations. Therefore, the research on Codazzi tensors is vast and has found many applications, and it would require a very long read to cover all aspects, even superficially. Some of these results can be found in the monograph [2] which was published in 1987. On the other hand, there are many papers on the geometry of Codazzi tensors [5] -  [10] which were published in subsequent years.

\noindent In turn, we introduced in [11] the notion of Codazzi ${\rm p}$-tensors $\left({\rm p}\ge 2\right)$ which extends the well known concept for $p=2$. Let us recall that a \textit{Codazzi ${\rm p}$-tensor} or, in other words, a \textit{higher order Codazzi tensor} is a $C^{\infty } $-section \textit{T} of the vector bundle ${\rm S}^{{\rm p}} M$on \textit{${\rm M}$ } satisfying the following equation:
\begin{equation} \label{GrindEQ__1_2_} 
\left(\nabla _{X_{0} } T\right)\, \left(X_{1} ,\, X_{2} ,\, ...,\, X_{p} \right)=\left(\nabla _{X_{1} } T\right)\, \left(X_{0} ,\, X_{2} ,\, ...,\, X_{p} \right) 
\end{equation} 
for any tangent vector fields $X_{0} ,\, X_{1} ,\, ...,\, X_{p} $ on\textit{ $M$}. The theory of higher order Codazzi tensors was developed in the papers from the following list [12] - [15].

\noindent In the present paper we study the question of nonexistence of higher order Codazzi tensors and classical Codazzi tensors on complete and compact Riemannian manifolds, in particular. The classical \textit{Bochner technique} [16]; [17]; [18, pp. 333-363]; [19] and its generalized version [20] will help us to accomplish this task. We must recall here that the classic Bochner technique is an analytical method to obtain vanishing theorems for some topological or geometrical invariants on a compact (without boundary) Riemannian manifold, under some curvature assumption. The proofs of such theorems apply the \textit{Bochner maximum principle} and the \textit{theorem of Green} [16, pp. 30-31]. We will also use a generalized version of the Bochner technique in the present paper. Therefore, we will use the \textit{maximum principle of Hopf} [21, p. 47], Yau, Li and Schoen results on the connections between the geometry of a complete smooth manifold and the global behavior of its subharmonic functions [22]; [23]; [24]. We have demonstrated in our papers [25]; [26] and [27] that this generalized version of the Bochner technique is effective for the differential geometry ``in the large''. In addition, the theorems and corollaries of the present paper will illustrate the effectiveness of this method for studying the global geometry of higher order Codazzi tensors and, in particular, Codazzi tensor of order 2 on complete and compact Riemannian manifolds (see our Theorem 1, Corollary 2 and Theorem 3). 

\noindent It is well known that the Ricci tensor $Ric$ is a Codazzi 2-tensor on an \textit{n}-dimensional $\left(n\ge 3\right)$ locally conformally flat manifold \textit{$\left({\rm M,g}\right)$ }with constant (not necessarily zero) scalar curvature ${\rm s}={\rm trace}_{{\rm g}} Ric$ [28]. By using this fact, we will give an application of the Bochner technique to the global geometry of complete locally conformally flat Riemannian manifolds (see our Corollary 5 and Theorem 4). On the other hand, if $\left({\rm M,}\, {\rm g}\right)$ is a minimal hypersurface in the standard sphere ${\rm S}^{{\rm n}+1} $, then its second fundamental form ${\rm S}$ is a traceless Codazzi 2-tensor [29, p. 388]. By using this theorem, we will give an application of the Bochner technique to the global geometry of complete minimal hypersurfaces in a sphere (see Theorem 5). 

\noindent The theorems and corollaries of this paper supplement our results from [1] and the results of other authors from [2, Theorem 16.9]; [13]; [28] and [30]. 

\noindent \textbf{2. Main results}

\noindent Everywhere in this paper we consider a higher order Codazzi tensor $T$ as a smooth section of the subbundle $S_{0}^{p} M$ of the bundle $S^{{\rm p}} M$ on a Riemannian manifold \textit{$\left({\rm M,g}\right)$ }defined by the condition
\[{\rm trace}_{{\rm g}} \, T=\sum _{i=1,...,n}T\, \left(\, e_{i} ,\, e_{i} ,\, X_{3} ,\, ...,\, X_{p} \right) =0\] 
for $T\in C^{\infty } S_{0}^{p} M$ and orthonormal basis $\left\{\, e_{i} \, \right\}$ of $T_{x}M $ at an arbitrary point $x\in M$. In this case, ${\rm T}$ is called \textit{traceless}. It can be proved that an arbitrary traceless Codazzi ${\rm p}$-tensor \textit{T} is a \textit{divergence-free }tensor field, i.e.\textit{ }$\delta \, T=0$ for the formal adjoint operator $\delta {\rm :}C^{\infty } \left(\otimes ^{p+1} {\rm T}^{*} M\right)\to C^{\infty } \left(\otimes ^{p} {\rm T}^{*} M\right)$ of the covariant derivative $\nabla {\rm :}C^{\infty } \left(\otimes ^{p} {\rm T}^{*} M\right)\to C^{\infty } \left(\otimes ^{p+1} {\rm T}^{*} M\right)$ where ${\rm T}^{*} M$ is the cotangent bundle on ${\rm M}$(see [2, p. 54]; [29]). 

\noindent \textbf{Remark 1}. In [12] the dimension of the vector space of traceless Codazzi ${\rm p}$-tensors $\left({\rm p}\ge 2\right)$ on compact Riemannian surfaces of genus $\gamma $ was determined. It depends only on ${\rm p}$ and $\gamma $. Additionally the result was extended to genus zero.~

The following \textit{Bochner-Weitzenb\"{o}ck formula} for an arbitrary Codazzi ${\rm p}$-tensor ${\rm T}\in {\rm C}^{\infty } S_{0}^{{\rm p}} M$ holds
\begin{equation} \label{GrindEQ__2_1_} 
\frac{1}{2} \, \Delta _{{\rm B}} \, \left\| \, {\rm T}\, \right\| ^{2} =Q_{P} \, \left(T{\rm ,}\, {\rm T}\right)+\left\| \, \nabla \, T\, \right\| ^{2} ,                                
\end{equation} 
where $\Delta _{B} ={\rm div}\circ \, {\rm grad}$ is the \textit{Laplace-Beltrami operator}, $\left\| \, {\rm T}\, \right\| ^{2} =g\, \left(\, T,\, T\right)$,  $\left\| \, \nabla {\rm T}\, \right\| ^{2} =g\, \left(\, \nabla T,\, \nabla T\right)$ and ${\rm Q}_{{\rm P}} $ is a quadratic form ${\rm Q}_{{\rm p}} {\rm :}\; S_{0}^{p} M\otimes S_{0}^{p} M\to $$\mathbb{R}$ which can be algebraically expressed through the curvature tensor $R$ and the Ricci tensor $Ric$ of \textit{$\left({\rm M,g}\right)$.} 

\noindent The curvature tensor ${\rm R}$ of \textit{$\left({\rm M,g}\right)$} induces an algebraic \textit{curvature operator} $\mathop{R}\limits^{\circ } :S_{0}^{2} M\to S_{0}^{2} M$. The symmetries of the curvature tensor $R$ imply that $\mathop{R}\limits^{\circ } $ is a selfadjoint operator, with respect to the point-wise inner product on $S_{0}^{2} M$. That is why, the eigenvalues of $\mathop{R}\limits^{\circ } $ are all real numbers at each point $x\in M$. Thus, we say $\mathop{R}\limits^{\circ } $ is \textit{positive semidefinite} (resp. \textit{strictly positive}), or simply $\mathop{R}\limits^{\circ } \ge 0$ (resp. $\mathop{R}\limits^{\circ } >0$), if all the eigenvalues of $\mathop{R}\limits^{\circ } $ are nonnegative (resp. strictly positive). In the next paragraph we will prove that $Q_{P} \left(T,\, T\right)\ge 0$ for an arbitrary ${\rm T}\in {\rm C}^{\infty } S_{0}^{p} M$ if the curvature operator $\mathop{R}\limits^{\circ } $ is positive semidefinite on $S_{0}^{2} M$ and $Q_{P} \left(T,\, T\right)>0$ for an arbitrary nonzero ${\rm T}\in {\rm C}^{\infty } S_{0}^{p} M$ if the curvature operator $\mathop{R}\limits^{\circ } $ is positive definite on $S_{0}^{2} M$. 

\noindent \textbf{Remark 2. }The curvature operator $\mathop{R}\limits^{\circ } $ has been studied in many papers and monographs. It is famous for its numerous applications [2, pp. 51-52; 346-347]; [31]; [32]; [34]; [35]. Beside the curvature operator $\mathop{R}\limits^{\circ } $ there is a curvature operator $\bar{{\rm R}}{\rm :}\Lambda ^{2} M\to \Lambda ^{2} M$ (see [18, pp. 83]). It is also widely used in Riemannian geometry. Examples are given by the well known Gallot-Meyer theorem [18, p. 351] on harmonic $p$-forms $\omega \in {\rm C}^{\infty } \Lambda ^{p} M$ on a compact manifold \textit{$\left({\rm M,g}\right)$ }and its $p$-th Betti number $\beta _{{\rm p}} \left(M\right)$, and theorems from [25] and [26]. Therefore $\mathop{R}\limits^{\circ } $ is also referred to as the \textit{curvature operator of the second kind} [34].

\noindent Taking into account \eqref{GrindEQ__2_1_} and applying the "Hopf maximum principle" [9, p. 47] or, in other words, ``strong maximal principal of Hopf'' we will prove the following statement.

\noindent \textbf{LEMMA 1}.\textit{ Let the curvature operator }$\mathop{R}\limits^{\circ } \ge 0$\textit{ at} \textit{any point of a connected open domain $U$ of $\left({\rm M,g}\right)$ and ${\rm T}$ be a Codazzi }$p$-\textit{tensor }($p\ge 2$) \textit{defined at any point of ${\rm U}$. If} \textit{the scalar function $\left\| \, {\rm T}\, \right\| ^{2} $ has a local maximum at some point of }${\rm U}$\textit{, then $\left\| \, {\rm T}\, \right\| ^{2} $ is a constant function and ${\rm T}$ is invariant under parallel translations} \textit{in }${\rm U}$\textit{. Moreover, if }$\mathop{R}\limits^{\circ } >0$\textit{ at some point }$x\in U$\textit{, then $T\equiv 0$.}

\noindent Let \textit{$\left({\rm M,g}\right)$ }be a compact Riemannian manifold with the curvature operator $\mathop{R}\limits^{\circ } \ge 0$. Then there exists a point ${\rm x}\in {\rm M}$ at which the function \textit{$\left\| \, {\rm T}\, \right\| ^{2} $ }attains the maximum. At the same time, \textit{$\left\| \, {\rm T}\, \right\| ^{2} $} satisfies the condition $\Delta \, \left\| \, {\rm T}\, \right\| ^{2} \ge 0$ everywhere in \textit{$\left({\rm M,g}\right)$}. If, moreover, $\mathop{R}\limits^{\circ } >0$ at some point $x\in {\rm M}$, then $T\equiv 0$ everywhere on \textit{$\left({\rm M,g}\right)$}. As a result, we obtain the corollary that was proved in our paper [11]. Moreover, it is a generalization of the Berger-Ebin theorem on an arbitrary Codazzi 2-tensor with constant trace on a compact Riemannian manifold [2, p. 436]; [29].

\noindent \textbf{COROLLARY 1.} \textit{Every traceless Codazzi $p$-tensor }(\textit{$p\ge 2$})\textit{ on a compact Riemannian manifold $\left({\rm M,g}\right)$ with nonnegative curvature operator is invariant under parallel translations. Moreover, if }$\mathop{R}\limits^{\circ } >0$\textit{ at some point }$x\in {\rm M}$\textit{, then there is no non-zero traceless Codazzi $p$-tensor }(\textit{$p\ge 2$})\textit{.}

\noindent \textbf{Remark 3}. In [37] was proved that \textit{g }($\mathop{R}\limits^{\circ } \, \left({\rm T}\right),\, T$) $\geq$ 0 for any $T\in C^{\infty } S_{0}^{p} M$ and ${\rm p}\ge 2$ if ${\rm sec}\ge 0$ for the sectional curvature ${\rm sec}$ of $\left(M,g\right)$. Therefore, we can reformulate our Lemma 1 and Corollary 1 by using this statement. In particular, we can state that every traceless higher order Codazzi tensor on a compact Riemannian manifold with nonnegative sectional curvature is invariant under parallel translations.

\noindent Let us formulate a theorem that supplements Theorem 2 from our paper [11] which was proved for a compact Riemannian manifold.

\noindent \textbf{THEOREM}\textit{ }\textbf{1}\textit{. Let $\left({\rm M,g}\right)$ be a complete noncompact Riemannian manifold with nonnegative curvature operator }$\mathop{R}\limits^{\circ } $\textit{. Then there is no non-zero traceless Codazzi $p$-tensor} ($p\ge 2$)\textit{ on $\left({\rm M,g}\right)$ such that }$\int _{{\rm M}}\left\| \, {\rm T}\, \right\| ^{{\rm q}} dVol_{g} <+\infty  $ \textit{for some $q\ge 1$.}

\noindent A Riemannian manifold \textit{$\left({\rm M,g}\right)$} is \textit{locally conformally flat} if, for any ${\rm x}\in {\rm M}$, there exists a neighborhood $U$of $x$ and $C^{\infty } $-function $f$ on $U$ such that $\left(\, U,\, e^{2f} g\right)$ is flat [2, p. 60]. We remind that a manifold \textit{$\left({\rm M,g}\right)$ }of dimension ${\rm n}\; \; \left({\rm n}\ge 4\right)$ is locally conformally flat if and only if its Weyl tensor ${\rm W}$ is identically zero [2, p. 60]. The Schouten tensor ${\rm Sch}$ plays an important role in the description of such manifolds [38]. This tensor has the form ${\rm S}ch=\left(n-2\right)^{-1} \, \left(\, R{\rm ic}-s\, \left(2n-2\right)^{-1} g\right)$ for the Ricci tensor ${\rm Ric}$ and the scalar curvature ${\rm s}={\rm trace}_{{\rm g}} \, {\rm Ric}$ of \textit{$\left({\rm M,g}\right)$}. Let us formulate the following statement.

\noindent \textbf{COROLLARY} \textbf{2}. \textit{Let $\left({\rm M,g}\right)$ be a complete noncompact locally conformally flat Riemannian manifold of dimension ${\rm n}\; \; \left({\rm n}\ge 4\right)$ with the nonnegative definite Schouten tensor}. \textit{Then there is no non-zero traceless Codazzi $p$-tensor} ($p\ge 2$)\textit{ on $\left({\rm M,g}\right)$ such that }$\int _{{\rm M}}\left\| \, {\rm T}\, \right\| ^{{\rm q}} dVol_{g} <+\infty  $ \textit{for some $q\ge 1$.}

\noindent \textbf{Remark 4}. Corollary 2 supplements to some results on compact locally conformally flat Riemannian manifolds from [13]. We note that the condition of the nonnegative definiteness of the Schouten tensor in the case of constant positive scalar curvature means for a complete manifold that it is compact [18, p. 251].

\noindent Let us now consider a non-zero Codazzi 2-tensor on a Riemannian manifold, i.e. a symmetric bilinear form ${\rm T}\in C^{\infty } S^{2} M$ satisfying the \textit{Codazzi equation} \eqref{GrindEQ__1_1_}. 

\noindent First, we consider a symmetric bilinear form ${\rm T}\in C^{\infty } S^{2} M$ as a 1-form with values in the cotangent bundle $T^{*} M$. This bundle comes equipped with the Levi-Civita covariant derivative $\nabla $, thus there is an induced exterior differential $d^{\nabla } {\rm :C}^{\infty } S^{2} M\to C^{\infty } \left(\Lambda ^{2} M\otimes T^{*} M\right)$ on $T^{*} M$-valued differential one-forms such as ${\rm d}^{\nabla } T\, \left(X,\, Y,\, Z\right)=\left(\nabla _{X} {\rm T}\right)\, \left(Y,\, Z\right)-\left(\nabla _{Y} T\right)\, \left(X,\, Z\right)$ for any tangent vector fields ${\rm X,}\, {\rm Y,}\, {\rm Z}$ on ${\rm M}$ and an arbitrary ${\rm T}\in C^{\infty } S^{2} M$ [2, p. 355]; [18, pp. 349-350]; [39]. In this case, a symmetric bilinear form ${\rm T}\in C^{\infty } S^{2} M$ is a Codazzi 2-tensor if and only if $d^{\nabla } {\rm T}$ vanishes [18, p. 350]; [39]. 

\noindent On the other hand, P. Peterson called [18, p. 350] a symmetric bilinear form ${\rm T}\in C^{\infty } S^{2} M$ \textit{harmonic} if $T\in {\rm Ker}\, d^{\nabla } \bigcap {\rm Ker}\, \delta $. In this case, ${\rm T}$ is a divergence free Codazzi 2-tensor. In addition, from \eqref{GrindEQ__1_1_} we obtain the equation $\delta \, T=-\, \; d\left(\, {\rm trace}_{g} T\right)$ for an arbitrary Codazzi 2-tensor ${\rm T}\in {\rm C}^{\infty } S^{2} M$. Therefore, we can conclude that a bilinear form ${\rm T}\in C^{\infty } S^{2} M$ is harmonic if and only if it is a Codazzi 2-tensor with constant trace [18, p. 350]. Moreover, the following statement holds.

\noindent \textbf{THEOREM 2.} \textit{The vector space of harmonic symmetric bilinear forms} \textit{on a compact Riemannian manifold is finite dimensional.}

\noindent \textbf{Remark 5. }J. P.\textbf{ }Bourguignon proved [40, p. 281] that a compact orientable Riemannian four-manifold admitting a non-trivial Codazzi tensor of order 2 with constant trace (harmonic symmetric bilinear form) must have signature zero. \textbf{}

\noindent Everywhere in the following we will consider a harmonic symmetric bilinear form or, in other words, we will consider a Codazzi 2-tensor with constant trace.

\noindent We have the Bochner-Weitzenb\"{o}ck formula 
\begin{equation} \label{GrindEQ__2_2_} 
\frac{1}{2} \, \Delta _{B} \, \left\| \, {\rm T}\, \right\| ^{2} =Q_{2} \, \left(T,\, T\right)+\left\| \, \nabla \, T\, \right\| ^{2} ,                                 
\end{equation} 
for an arbitrary harmonic symmetric bilinear form. Here the sign of the quadratic form  $Q_{2} \left(T,\, T\right)$ depends on the sign of the sectional curvature ${\rm sec}$ of $\left({\rm M,g}\right)$ [29]. We remind here that an arbitrary Codazzi 2-tensor \textit{T} on \textit{$\left({\rm M,g}\right)$} commutes with the Ricci tensor $Ric$ at each point ${\rm x}\in {\rm M}$ [1]; [2, p. 439]. Therefore, the eigenvectors of an arbitrary Codazzi tensor $T$ determine the principal directions of the Ricci tensor at each point ${\rm x}\in {\rm M}$ [41, pp.113-114]. The converse is also true. Then taking into account of \eqref{GrindEQ__2_2_} and using the "Hopf maximum principle", we will prove in the next paragraph that the following lemma holds.

\noindent \textbf{LEMMA 2}. \textit{Let $U$be a connected open domain $U$ of a Riemannian manifold $\left({\rm M,g}\right)$ and ${\rm T}$ be a harmonic symmetric bilinear form defined at any point of ${\rm U}$. If the sectional curvature ${\rm sec}\, \left(\, {\rm e}_{{\rm i}} \wedge \, e_{j} \right)\ge 0$ for all vectors of the orthonormal basis $\left\{\, e_{i} \right\}$ of $T_{x} M$ which is determined by the principal directions of the Ricci tensor ${\rm Ric}$ at an arbitrary point $x\in U$ and $\left\| \, {\rm T}\, \right\| ^{2} $ has a local maximum in the domain U, then $\left\| \, {\rm T}\, \right\| ^{2} $ is a constant function and ${\rm T}$ is invariant under parallel translations} \textit{in }${\rm U}$\textit{. Moreover, if ${\rm sec}\, \left(\, {\rm e}_{{\rm i}} \wedge \, e_{j} \right)>0$ at some point }$x\in U$\textit{, then ${\rm T}$ is trivial.}

\noindent If \textit{$\left({\rm M,g}\right)$} is a compact manifold and a harmonic symmetric bilinear form \textit{T} is given in a global way on \textit{$\left({\rm M,g}\right)$} then due to the "Bochner maximum principle" for compact manifold it follows the classical Berger-Ebin theorem [2, p. 436] and [29] which is a corollary of our Lemma 2 (and see also Remark 3).

\noindent \textbf{COROLLARY 3}. \textit{Every harmonic symmetric bilinear form }${\rm T}\in {\rm C}^{\infty } S^{2} M$\textit{ on a compact Riemannian manifold $\left({\rm M,g}\right)$ with nonnegative sectional curvature is invariant under parallel translations. Moreover, if ${\rm sec}>0$ at some point, then }${\rm T}\in {\rm C}^{\infty } S^{2} M$\textit{ is trivial.}

\noindent \textbf{Remark 6}. It is well known that every parallel symmetric tensor field $T\in C^{\infty } S^{2} M$ on a connected locally irreducible Riemannian manifold \textit{$\left({\rm M,g}\right)$ }is proportional to $g$, i.e. \textit{$T=\lambda \, g$ }for some constant $\lambda $. Due to this the second parts of Corollary\textit{ }3 can be reformulated in the following form: Moreover, if \textit{$\left({\rm M,g}\right)$} a connected locally irreducible Riemannian then $T$ is trivial.

\noindent For example, let $\left({\rm M,}\; {\rm g}\right)$ be a \textit{Riemannian symmetric space of compact type} that is a compact Riemannian manifold with non-negative sectional curvature and positive-definite Ricci tensor (see [42, p. 256]). Moreover, if a Riemann symmetric space of compact type is a locally irreducibility Riemannian manifold \textit{$\left({\rm M,g}\right)$ }then it is a compact Riemannian manifold with positive sectional curvature [43]. Therefore, we can formulate the following corollary.

\noindent \textbf{COROLLARY\textit{ }4}\textit{.} \textit{Every harmonic symmetric bilinear form} \textit{on a Riemannian symmetric manifold of compact type is invariant under parallel translations. If in addition to the above mentioned the manifold is locally irreducible then harmonic symmetric bilinear forms are trivial.}

\noindent The following theorem supplements the Berger-Ebin theorem [2, p. 436] and [29] for the case of a complete noncompact Riemannian manifold.

\noindent \textbf{THEOREM} \textbf{3}. \textit{Let $\left({\rm M,g}\right)$ be a connected complete noncompact Riemannian manifold with nonnegative sectional curvature. Then there is no a non-zero harmonic symmetric bilinear form ${\rm T}$} \textit{which satisfies the condition }$\int _{{\rm M}}\left\| \, T\, \right\| ^{{\rm q}} \, dVol_{g} <+\infty  $\textit{ at least for one} $q\ge 1$.

\noindent \textbf{Remark 7}. In the case of a locally conformally flat \textit{n}-dimensional $\left(\, n\ge 4\, \right)$ Riemannian manifold \textit{$\left({\rm M,g}\right)$} the following equalities hold [28]; [44]
\[\sec \, \left(\, e_{i} \wedge \, e_{j} \right)=\frac{1}{n-2} \, \left(r_{i} +r_{j} -\frac{1}{n-1} \left(\, r_{\, 1} +...+r_{\, n} \right)\right);    \sec \, \left(e_{i} \wedge \, e_{j} \right)=\frac{1}{n-2} \, \left(\, \lambda _{i} +\lambda _{j} \right)\] 
where \textit{$\left\{\, e_{i} \right\}$} is a orthonormal basis of $T_{x} M$ at an arbitrary point $x\in M$ such  that  \textit{$Ric\, \left(e_{i} ,\, e_{j} \right)=r_{i} \, \delta _{ij} $},\textit{ $Sch\, \left(e_{i} ,\, e_{j} \right)=\lambda _{i} \, \delta _{ij} $} and $\lambda _{{\rm i}} =\left(n-2\right)^{-1} \left(\, r_{i} -\left(2n-2\right)^{-1} s\, \right)$. Due to these equations we can formulate analogues of the Lemma 2 and the Corollary\textit{ }3 for a locally conformally flat Riemannian manifold where the inequality \textit{$\sec \, \left(\, e_{i} \wedge \, e_{j} \right)\ge 0$ }can be replaced by \textit{$\lambda _{i} +\lambda _{j} \ge 0$ }or \textit{$r_{i} +r_{j} \ge \left({\rm n}-1\right)^{-1} \left(\, r_{\, 1} +...+r_{\, n} \right)$ } for any ${\rm i}\ne {\rm j}$.

\noindent For an \textit{n}-dimensional $\left(n\ge 3\right)$ locally conformally flat manifold \textit{$\left({\rm M,g}\right)$ }with constant (not necessarily zero) scalar curvature its Ricci tensor ${\rm Ric}$ is a Codazzi 2-tensor with constant trace or, in other words, a harmonic symmetric bilinear form [2, p. 444]; [28]. Therefore, from Theorem 3 we conclude that the following corollary holds.

\noindent \textbf{COROLLARY 5\textit{.}} \textit{Let $\left({\rm M,g}\right)$ be an $n$-dimensional $\left(\, n\ge 3\, \right)$ connected complete noncompact locally conformally flat Riemannian manifold with nonnegative sectional curvature and constant scalar curvature. If $\left({\rm M,g}\right)$ is not locally flat then }$\int _{{\rm M}}\left\| \, {\rm Ric}\, \right\| ^{\, q} \, dVol_{g} =+\infty  $\textit{ for the Ricci tensor ${\rm Ric}$ and an arbitrary }$q\ge 1$\textit{.}

\noindent \textbf{Remark 8}. Our Corollary 5\textit{ }is a supplement to the theorem of M. Tani [28]. 

\noindent We strengthen the Corollary\textit{ }5 by proving the validity of the following theorem.

\noindent \textbf{THEOREM 4}.\textit{ Let $\left({\rm M,g}\right)$ be an $n$-dimensional $\left(\, n\ge 3\, \right)$ connected complete locally conformally flat Riemannian manifold such that }$\, \left\| \, {\rm Ric}\, \right\| ^{2} <\left(n-1\right)^{-1} s^{2} $ \textit{for} \textit{the Ricci tensor ${\rm Ric}$ and the positive constant scalar curvature ${\rm s}$}. \textit{If one of the following conditions is satisfied}:\textit{ }

\noindent (i) \textit{$\, \left\| \, {\rm Ric}\, \right\| ^{2} $ has a global maximum point};

\noindent (ii) $\int _{{\rm M}}\left\| \, {\rm Ric}\, \right\| ^{{\rm q}} \, dVol_{g} <+\infty  $\textit{ at least for one} $q\ge 2$;

\noindent (iii)\textit{ $\left({\rm M,g}\right)$ is a parabolic manifold, }

\noindent \textit{then $\left({\rm M,g}\right)$ is a spherical space form.}

\noindent \textbf{Remark 9. }Theorem 4 is a supplement to the following S. Goldberg theorem [30]. 

Let \textit{$\left({\rm M,}\, {\rm g}\right)$} be an ${\rm n}$-dimensional connected complete minimal hypersurface in the standard sphere $\left(\, {\rm S}^{{\rm n}+1} ,\, g_{0} \right)$ where ${\rm S}^{{\rm n}+1} \subset $ $\mathbb{R}$${}^{B}$ for ${\rm m}>{\rm n}+1$. In this case, we know that the second fundamental form ${\it S}$ of the hypersurface \textit{$\left({\rm M,}\, {\rm g}\right)$ }is a traceless Codazzi tensor [29, p. 388]. Therefore, \eqref{GrindEQ__2_2_} can be rewritten in the form [45]
\begin{equation} \label{GrindEQ__2_3_} 
\frac{1}{2} \, \Delta _{B} \, \left\| \, {\it S}\, \right\| ^{2} =\left\| \, {\rm S}\, \right\| ^{2} \left(\, n-\left\| \, {\rm S}\, \right\| ^{2} \right)+\left\| \, \nabla \, S\, \right\| ^{2}  
\end{equation} 
The following theorem holds. 

\noindent \textbf{Theorem 5.} \textit{Let $\left({\rm M,}\, {\rm g}\right)$ be a connected complete minimal hypersurface in the standard sphere} $\left(\, {\rm S}^{{\rm n}+1} ,\, g_{0} \right)$ \textit{such that its second fundamental form ${\rm S}$ satisfies the inequality $\left\| \, {\rm S}\, \right\| ^{2} \le {\rm n}$.} \textit{If one of the following conditions is satisfied}:\textit{ }

\noindent  (i) \textit{$\, \left\| \, {\rm S}\, \right\| ^{2} $ has a global maximum point};

\noindent (ii) \textit{$\int _{{\rm M}}\left\| \, {\rm S}\, \right\| ^{{\rm q}} \, dVol_{g} <+\infty  $ at least for one $q\ge 2$;}

\noindent (iii)\textit{ $\left({\rm M,g}\right)$ is a parabolic manifold }

\noindent \textit{then one of the following occurs}: \textit{}

\noindent (i) \textit{$\left({\rm M,}\, {\rm g}\right)$ is an equator of }$\left(\, {\rm S}^{{\rm n}+1} ,\, g_{0} \right)$;

\noindent (ii) \textit{$\left({\rm M,}\, {\rm g}\right)$is locally isometric to a generalized Clifford torus.  }

\noindent \textbf{Remark 10. }If \textit{$\left({\rm M,}\, {\rm g}\right)$} is a complete Riemannian manifold with finite volume then the conditions (ii) of Theorem 4 and Theorem 5 are satisfied. For example, from \textit{$\left\| \, {\rm S}\, \right\| ^{2} \le {\rm n}$} we obtain \textit{$\int _{{\rm M}}\left\| \, S\, \right\| ^{\, 2} dVol_{g} \le n\, \int _{{\rm M}}\, dVol_{g} \le n\, Vol_{g} \left(M\right)<+\infty   $.}

\pagebreak
\noindent \textbf{3. Proofs of the statements}

\noindent Let us deduce the Weitzenb\"{o}ck formula \eqref{GrindEQ__2_1_} for a traceless Codazzi ${\rm p}$-tensor \textit{T} for $p\ge 2$. For this purpose, we remind that the local components of the Ricci tensor $Ric$ of the manifold $\left(M, g\right)$ are calculated from the identity $Ric\, \left(\partial _{i} ,\partial _{j} \right)=R_{ij} =R_{ikj}^{k} $ for the local components $R_{ikj}^{l} $ of the curvature tensor \textit{R} determined from the equality $R_{lij}^{k} X^{l} =\nabla _{i} \nabla _{j} X^{k} -\nabla _{j} \nabla _{i} X^{k} $ where $\nabla _{i} =\nabla _{\partial _{{\rm i}} } $ and $X=X^{i} \partial _{i} $. Let us denote by $s:=g^{ij} R_{ij} $ the scalar curvature of the metric \textit{g }for $\left(g^{ij} \right)=\left(g_{ij} \right)\, ^{-1} $. Then direct calculations give 
\[\frac{1}{2} \, \Delta _{B} \, \left\| \, T\, \right\| ^{2} {\rm :}=\frac{1}{2} g^{kl} \, \nabla _{k} \, \nabla _{l} \, \left(T^{{\rm k}_{1} k_{2} ...k_{p} } T_{k_{1} k_{2} ...k_{p} }^{} \right)=\] 
\[=\left(g^{kl} \, \nabla _{k} \, \nabla _{l} \, T_{k_{1} k_{2} ...k_{p} }^{} \right)\, \, T^{{\rm k}_{1} k_{2} ...k_{p} } +g^{kl} \left(\nabla _{k} \, T_{k_{1} k_{2} ...k_{p} }^{} \right)\, \left(\nabla _{l} \, \, T^{{\rm k}_{1} k_{2} ...k_{p} } \right)=\] 
\[=R_{ij} \, T^{ik_{2} ...k_{p} } T_{k_{2} ...k_{p} }^{j} -\left({\rm p}-1\right)\, R_{ijkl} T^{ikk_{3} ...k_{p} } T^{jl} _{k_{3} ...k_{p} } +g^{kl} \left(\nabla _{k} \, T_{k_{1} k_{2} ...k_{p} }^{} \right)\, \left(\nabla _{l} \, \, T^{{\rm k}_{1} k_{2} ...k_{p} } \right)=\] 
\[=Q_{p} \, \left(T{\rm ,}\, {\rm T}\right)+\left\| \, \nabla \, T\, \right\| ^{2} \] 
where the quadratic form $Q_{p} \, \left(T{\rm ,}\, {\rm T}\right)$ has the form 
\begin{equation} \label{GrindEQ__3_1_} 
{\rm Q}_{p} \left(\, T\, ,T\right)=R_{ij} \, T^{ik_{2} ...k_{p} } T_{k_{2} ...k_{p} }^{j} -\left({\rm p}-1\right)\, R_{ijkl} T^{ikk_{3} ...k_{p} } T^{jl} _{k_{3} ...k_{p} }  
\end{equation} 
for components $T_{k_{1} ...k_{p} } $ of an arbitrary Codazzi ${\rm p}$-tensor ${\rm T}\in {\rm C}^{\infty } S_{0}^{p} M$with respect to a local coordinate system $x^{1} ,\, ...,\, x^{n} $. In deriving the formula we have used the condition of divergence-free of the Codazzi ${\rm p}$-tensor ${\rm T}$ and well-known Ricci identity (11.16) from [41].

\noindent The curvature operator $\mathop{R}\limits^{\circ } :S_{0}^{2} M\to S_{0}^{2} M$ is defined by equations [2, pp. 51-52] 
\begin{equation} \label{GrindEQ__3_2_} 
\mathop{R}\limits^{\circ } \left(\varphi _{il} \right)=R_{ijkl} \, {\rm T}^{jk} =g^{km} g^{jp} R_{ijkl} T_{mp}  
\end{equation} 
for local components $T_{ij} \; $ of arbitrary $T\in $$S_{0}^{2} M$. Everywhere else we assume that the curvature operator $\mathop{R}\limits^{\circ } $ is nonnegative definite on an arbitrary section of the bundle $S_{0}^{2} M$, i.e. the inequality $R_{ijkl} \, T^{il} T_{\quad }^{jk} \ge 0$ is true for an arbitrary ${\rm T}\in {\rm S}_{0}^{2} M $and then  $R_{ijkl} \, T^{il\, k_{3} ...k_{p} } T_{\quad k_{3} ...k_{p} }^{jk} \ge 0$ for any ${\rm T}\in {\rm S}_{0}^{p} M$. As a result, the second term in \eqref{GrindEQ__3_1_} will be nonpositive, i.e. $\left({\rm p}-1\right)\, R_{ijkl} T^{ikk_{3} ...k_{p} } T^{jl} _{k_{3} ...k_{p} } \le 0$. 

\noindent At an arbitrary point $ x\in M$ we choose orthogonal unit vectors $X,Y\in T_{x} M$ and define the tensor $\theta \in S_{0}^{2} M$ by the equality 
\[\theta =2^{-1} \left(X\otimes Y+Y\otimes X\right)\] 
then 
\[g\left(\mathop{R}\limits^{\circ } \left(\theta \right),\, \theta \right)=2g\left(R\left(X,Y\right)\, Y,\, X\right) = 2 \sec \, \left(X\wedge Y\right).\] 
That is why the sectional curvature of a manifold $\left(M, g\right)$ is everywhere nonnegative if the operator $\mathop{R}\limits^{\circ } $ is nonnegative definite on any section of the bundle $S_{0}^{2} M$. Let $X\in T_{x} M$ is a unit vector and we complete it to an orthonormal basis ${\rm X,}\, e_{2} ,\, ...\, ,e_{n} $ for $T_{x}M $ at an arbitrary point $x\in M$, then [18, p. 86]
\[Ric\left(X,X\right)=\sum _{a=2}^{n}sec\left(X\wedge \, e_{a} \right)\, . \] 
Therefore, the Ricci curvature is also nonnegative definite. Thus, if the operator $\mathop{R}\limits^{\circ } $ is nonnegative definite on sections of the bundle $S_{0}^{2} M$, then $Q_{p} \left(T,\, T\right)\ge 0$. As a result of \eqref{GrindEQ__2_2_} it follows that $\Delta _{{\rm B}} \, \left\| \, {\rm T}\, \right\| ^{2} \ge 0$, i.e. $\left\| \, {\rm T}\, \right\| ^{2} $ is a nonnegative subharmonic function. Moreover, we note that $Q_{p} \left(T_{x} ,\, T_{x} \right)>0$ for a nonzero ${\rm T}_{{\rm x}} \in S_{0}^{p} \left(T_{x} M\right)$ at some point $x\in M$ if $\mathop{R}\limits^{\circ } >0$ for all non-zero $\theta _{{\rm x}} \in S_{0}^{2} \left(T_{x} M\right)$ at this point.

\noindent Let us prove our Lemma 1. Consider a traceless Codazzi $p$-tensor $(p\ge 2)$ in the connected open domain $U$ of $\left({\rm M,g}\right)$ where the curvature operator $\mathop{R}\limits^{\circ } \ge 0$ then $Q_{p} \left(T,\, T\right)\ge 0$. And according to \eqref{GrindEQ__2_2_} we conclude that $\Delta _{{\rm B}} \, \left\| \, {\rm T}\, \right\| ^{2} \ge 0$, i.e. $\left\| \, {\rm T}\, \right\| ^{2} $ is a subharmonic function in the domain $U$. Suppose also that $\left\| \, {\rm T}\, \right\| ^{2} $ has a local maximum at some point ${\rm C}\in {\rm U}$ then according to the "Hopf's maximum principle" [9, p. 47] we have that $\left\| \, {\rm T}\, \right\| ^{2} $ is a constant function in the domain ${\rm U}\subset {\rm M}$. In this case, $\Delta \, \left\| \, T\, \right\| ^{2} =0$ and as a consequence of \eqref{GrindEQ__2_2_} we obtain $\left\| \, \nabla \, T\, \right\| ^{2} =0$ which means that the Codazzi $p$-tensor \textit{T }is parallel.

\noindent Let $\left\| \, T\, \right\| ^{2} ={\rm C}$ (where \textit{C} is a constant), then $T$ does not become zero anywhere in the domain $U$ and at the same time $\Delta _{{\rm B}} \, \left\| \, {\rm T}\, \right\| ^{2} =0$. If there is a point ${\rm x}\in {\rm U}$ where $\mathop{R}\limits^{\circ } >0$ then, as we stated above, $Q_{p} \left(T_{x} ,\, T_{x} \right)>0$. In this case, the left site of \eqref{GrindEQ__2_2_} is equal to zero and its right side is greater than zero, so it follows that ${\rm T}=0$. The Lemma 1 is proved.

\noindent Corollary 1 of our Lemma 1 does not require any proof. 

\noindent For the proof of our Theorem 1 in the case ${\rm q}=1$ we use the Theorem 1 from [23] according to which any nonnegative subharmonic function ${\rm f}\in {\rm C}^{\infty } {\rm M}$ on a connected complete noncompact Riemannian manifold with nonnegative sectional curvature must satisfy the condition $\int _{{\rm M}}f\, dVol_{g} =+\infty  $.

\noindent Let us suppose that $\left({\rm M,g}\right)$ is a complete noncompact Riemannian manifold with nonnegative curvature operator $\mathop{{\rm R}}\limits^{\circ } $ (and that is why with nonnegative sectional curvature) and with a globally defined non-zero Codazzi ${\rm p}$-tensor  ${\rm T}\in {\rm C}^{\infty } S_{0}^{p} M$ for an arbitrary $p\ge 2$. 

\noindent By direct calculation we find the following
\[\frac{1}{2} \, \Delta _{B} \, \left\| \, {\rm T}\, \right\| ^{2} =\left\| \, {\rm T}\, \right\| \, \, \Delta _{B} \, \left\| \, {\rm T}\, \right\| +\left\| \, {\rm d}\, \left\| \, {\rm T}\, \right\| \, \right\| ^{2} .\] 
Then the equation \eqref{GrindEQ__2_1_} can be rewritten in the form
\[\left\| \, {\rm T}\, \right\| \, \, \Delta _{B} \, \left\| \, {\rm T}\, \right\| =Q_{p} \left(T,\, T\right)+\left\| \, \nabla \, T\, \right\| ^{2} -\left\| \, d\, \left\| \, {\rm T}\, \right\| \, \right\| ^{2} .\] 
By using the first Kato inequality $\left\| \, \nabla \, T\, \right\| ^{2} \ge \left\| \, d\, \left\| \, {\rm T}\, \right\| \, \right\| ^{2} $(see [46]), we can write 
\[\left\| \, {\rm T}\, \right\| \, \, \Delta _{B} \, \left\| \, {\rm T}\, \right\| \ge Q_{p} \left(T,\, T\right).\] 
Therefore, if ${\rm Q}_{{\rm p}} \left(T,\, T\right)\ge 0$ then we have $\Delta \, \left\| \, {\rm T}\, \right\| \ge 0$ and as a result of this $\, \left\| \, {\rm T}\, \right\| $ is a nonnegative subharmonic function. Now according to Theorem 1 from [23] we come to the conclusion that $\int _{{\rm M}}\left\| \, T\, \right\| \, dVol_{g} =+\infty  $ for an arbitrary non-zero Codazzi $p$-tensor ${\rm T}\in {\rm C}^{\infty } S_{0}^{p} M$. Thus, on a complete noncompact Riemannian manifold with nonnegative curvature operator (and that is why with nonnegative sectional curvature) there exists no non-zero Codazzi tensor ${\rm T}\in {\rm C}^{\infty } S_{0}^{p} M$ such that $\int _{{\rm M}}\left\| \, T\, \right\| \, dVol_{g} <+\infty  $.

\noindent For the proof of our Theorem 1 in the case ${\rm q}>1$ we use the Theorem 7 from [22]. According to this theorem any nonnegative subharmonic function on a connected complete noncompact Riemannian manifold satisfies the condition $\int _{{\rm M}}f^{{\rm q}} dVol_{g} =+\infty  $ for ${\rm q}>1$ or ${\rm f}={\rm const}$.

\noindent Let us suppose that a non-zero traceless Codazzi $p$-tensor ${\rm T}$ is globally defined on a complete noncompact Riemannian manifold \textit{$\left({\rm M,g}\right)$ }with nonnegative curvature operator $\mathop{{\rm R}}\limits^{\circ } $.  In this case$\, \left\| \, {\rm T}\, \right\| $ is a nonnegative subharmonic function. If in addition to the above mentioned ${\rm T}$ satisfies the condition $\int _{{\rm M}}\left\| \, {\rm T}\, \right\| ^{q} dVol_{g} <+\infty  $ for some ${\rm q}>1$, then due to the Yau theorem we conclude that $\left\| \, {\rm T}\, \right\| =C$ for some nonnegative constant ${\rm C}$. Therefore, the condition $\int _{{\rm M}}\left\| \, {\rm T}\, \right\| ^{q} dVol_{g} <+\infty  $ can be rewritten in the form $C^{q} \int _{{\rm M}}dVol_{g} <+\infty  $. On the other hand, every complete noncompact Riemannian manifold with nonnegative Ricci curvature or with nonnegative sectional curvature has infinite volume [22]; [23]. Therefore, from the inequality $C^{q} \int _{{\rm M}}dVol_{g} <+\infty  $ we conclude that $\left\| \, {\rm T}\, \right\| =C=0$. This completes our proof.

\noindent Let us prove Corollary 2. Suppose that \textit{$\left({\rm M,g}\right)$ }is a ${\rm n}$-dimensional $\left(\, n\ge 4\, \right)$ complete noncompact locally conformally flat Riemannian manifold. In this case, the curvature tensor ${\rm R}$ has the following local components [2, pp. 60-61]; [28]:
\begin{equation} \label{GrindEQ__3_3_} 
{\rm R}_{{\rm ijkl}} =\frac{1}{n\, -\, 2} \, \, \left(R_{jl} g_{ik} -R_{jk} g_{il} +R_{ik} g_{jl} -R_{il} g_{jk} \right)\, \, -\frac{s}{\left(n-1\right)\left(n-2\right)} \, \left(g_{jl} g_{ik} -g_{jk} g_{il} \right).   
\end{equation} 
Then from \eqref{GrindEQ__3_3_} we obtain the following equalities 
\[{\rm R}_{{\rm ijkl}} \, \theta ^{jk} \theta ^{il} =\frac{2}{n-2} \, \left(R_{ij} \, \theta ^{ik} \theta _{k}^{j} -\frac{s}{2\, \left(n-1\right)} \, \left\| \, \theta \, \right\| ^{2} \right)=\] 
\[=\frac{1}{n\, -\, 2} \, \left(\left(R_{jl} -\frac{s}{2\, \left(n-1\right)} \, g_{j{\rm l}} \right)\, g_{ik} +\left(R_{{\rm ik}} -\frac{s}{2\, \left(n-1\right)} \, g_{ik} \right)\, g_{jl} \right)\, \theta ^{il} \theta ^{jk} =\] 
\[=\, \left(S_{jl} \, g_{ik} +S_{ik} \, g_{jl} \right)\, \theta ^{il} \theta ^{jk} =2\, {\rm S}_{{\rm ij}} \theta ^{ik} \theta _{k}^{j} .\] 
for the local component ${\rm S}_{{\rm ik}} $ of the Schouten tensor $SA{\rm h}\in C^{\infty } S^{2} M$ and the local components $\theta _{ik} $ of an arbitrary traceless tensor $\theta \in S_{0}^{2} M$. Then from the condition of nonnegative definiteness of the Schouten tensor $Sch \in C^{\infty } S^{2} M$ we obtain the inequality $\mathop{{\rm R}}\limits^{\circ } \ge 0$. After that we should repeat the proof of Theorem 1.

\noindent Let us prove Theorem 2. We define the following differential operator: $\delta ^{*} {\rm :}C^{\infty } T^{*} M\to C^{\infty } S^{2} M$ by the equality $\left(\delta ^{*} \xi \right)\, \left(\, {\rm X,}\, {\rm Y}\right)=2^{-1} \left(\, \left(\nabla _{X} \xi \right)\, Y+\left(\nabla _{Y} \xi \right)\, X\right)$ for any tangent vector fields ${\rm X,}\, {\rm Y,}\, {\rm Z}$ on ${\rm M}$ and an arbitrary one-form $\xi \in {\rm C}^{\infty } T^{*} M$. Moreover, we denote by $\delta ^{\nabla } $ the adjoint operator of ${\rm d}^{\nabla } $ [29, p. 380; 388]; [2, p. 355]. If $\left({\rm M,}\, {\rm g}\right)$ is a compact manifold then the condition $T\in {\rm Ker}\, d^{\nabla } \bigcap {\rm Ker}\, \delta $ is equivalent to $T\in {\rm Ker}\, \Psi $ for  $\Psi =\delta ^{\nabla } d^{\nabla } +2\, \delta ^{*} \delta $ [29, p. 388]. The differential operator $\Psi $ is an elliptic operator with injective symbol and ${\rm Ker}\, \Psi ={\rm Ker}\, d^{\nabla } \bigcap {\rm Ker}\, \delta $ [29, p. 388]. Therefore, ${\rm Ker}\, \Psi $ is a finite-dimensional vector space of harmonic symmetric bilinear form on a compact manifold $\left({\rm M,}\, {\rm g}\right)$ [2, p. 464]. This completes our proof.

\noindent Let us prove Lemma 2. First, we rewrite the Bochner-Weitzenb\"{o}ck formula \eqref{GrindEQ__2_2_} in the form
\begin{equation} \label{GrindEQ__3_4_} 
\frac{1}{2} \, \Delta _{B} \, \left\| \, \bar{{\rm T}}\, \right\| ^{2} =Q_{2} \left(\bar{T},\, \bar{T}\right)+\left\| \, \nabla \, \bar{T}\, \right\| ^{2}  
\end{equation} 
where ${\rm Q}_{2} \left(\, \bar{T},\, \bar{T}\right)=R_{ij} \bar{T}^{ik} \bar{T}_{k}^{j} -R_{ijkl} \bar{T}^{ik} \bar{T}^{jl} $ for a traceless Codazzi 2-tensor $\bar{T}$. 

\noindent Second, if $T$ is a Codazzi tensor with constant trace, then $\bar{T}=T-{1 \mathord{\left/{\vphantom{1 n}}\right.\kern-\nulldelimiterspace} n} \, \, \left(trace_{g} \, T\right)\, g$ is a traceless Codazzi tensor such that $\Delta _{B} \, \left\| \, \overline{{\rm T}}\, \right\| ^{2} =\Delta _{B} \, \left\| \, T\, \right\| ^{2} $ and $\nabla \bar{{\rm T}}=\nabla T$. In addition, one can prove that  $Q_{2} \left(\, \bar{T},\, \bar{T}\right)=Q_{2} \left(\, T,\, T\right)=R_{ij} T^{ik} T_{k}^{j} -R_{ijkl} T^{ik} T^{jl} $. Then \eqref{GrindEQ__3_4_} can be rewritten in the form
\begin{equation} \label{GrindEQ__3_5_} 
\frac{1}{2} \, \Delta _{B} \, \left\| \, {\rm T}\, \right\| ^{2} =Q_{2} \left(T,\, T\right)+\left\| \, \nabla \, T\, \right\| ^{2}  
\end{equation} 
According to [2, p. 436] and [29, p. 388] we conclude that
\begin{equation} \label{GrindEQ__3_6_} 
{\rm Q}_{2} \left(\, T,\, T\right)=\sum _{{\rm i}\, <{\rm j}}{\rm sec} \left(e_{i} \wedge \, e_{j} \right)\, \left(\lambda _{i} -\lambda _{j} \right)\, ^{2} ,                                   
\end{equation} 
where \textit{$\left\{\, e_{i} \right\}$} is a such orthonormal basis of the tangent space $T_{x} M$ at an arbitrary point $x\in M$ that the Codazzi tensor $T_{x} \left(\, e_{i} ,\, e_{j} \right)=\lambda _{i} \left(x\right)\, \delta _{ij} $ where $\delta _{ij} $ is the Kronecker symbol and ${\rm sec}\, \left(\, {\rm e}_{{\rm i}} \wedge \, e_{j} \right)$ is the sectional curvature in the direction of subspace $\pi \left(x\right)\subset T_{x} M$ such that $\pi \left(x\right)=span\, \left\{\, e_{i} ,\, e_{j} \right\}$. It is known from [2, p. 439] and [1] that an arbitrary Codazzi tensor \textit{T} on a manifold \textit{$\left({\rm M,g}\right)$} commutes with its Ricci tensor $Ric$, and therefore the eigenvectors \textit{$\left\{\, e_{i} \right\}$} of the Codazzi tensor $T$ determine  the \textit{principal directions of the Ricci tensor} at each point ${\rm x}\in {\rm U}$[41, pp.113-114]. The converse is also true. Next, let us suppose that ${\rm sec}\, \left(\, {\rm e}_{{\rm i}} \wedge \, e_{j} \right)\ge 0$ in some connected open domain ${\rm U}\subset {\rm M}$ then ${\rm Q}\left(\, T\, \right)\ge 0$. Moreover, if there is a non-zero Codazzi tensor ${\rm T}$given in ${\rm U}\subset {\rm M}$ then from \eqref{GrindEQ__3_5_} we conclude that $\Delta _{B} \, \left\| \, T\, \right\| ^{2} \ge 0$, i.e. $\left\| \, T\right\| ^{2} $ is a nonnegative subharmonic function in $U$. Let us suppose $\left\| \, T\right\| ^{2} $ has a local maximum at some point ${\rm x}\in {\rm U}$ then $\left\| \, T\, \right\| ^{2} $ is a constant function in ${\rm U}\subset {\rm M}$ according to the "Hopf's maximum principle" [21, p. 47]. In this case, $\Delta _{B} \, \left\| \, T\, \right\| ^{2} =0$ and $\left\| \, \nabla \, T\, \right\| ^{2} =0$. Then the latter equation means that the Codazzi tensor ${\rm T}$ is a parallel tensor field.

\noindent Let $\left\| \, T\, \right\| ^{2} ={\rm C}$ (\textit{C} is a constant) then from \eqref{GrindEQ__3_5_} we obtain that $Q_{2} \left(T,\, T\right)+\left\| \, \nabla \, T\, \right\| ^{2} =0$. Since ${\rm sec}\, \left(\, {\rm e}_{{\rm i}} \wedge \, e_{j} \right)\ge 0$ it means that ${\rm Q}_{2} \left(\, T,\, T\right)=0$ and $\nabla \, T=0$. If there is a point ${\rm x}\in {\rm U}$ such that ${\rm sec}\, \left(\, {\rm e}_{{\rm i}} \wedge \, e_{j} \right)>0$ then from \eqref{GrindEQ__3_6_} we come to the conclusion that $\lambda _{1} ={\rm ...}=\lambda _{{\rm n}} =\lambda $ which is equivalent to $T={1 \mathord{\left/{\vphantom{1 n}}\right.\kern-\nulldelimiterspace} n} \, \, \left(trace\, T\right)\, g$ (see [2, p. 436]). Lemma 2 is proved.

\noindent Assume that the manifold \textit{$\left({\rm M,g}\right)$ }is compact and the Codazzi tensor is globally defined on \textit{$\left({\rm M,g}\right)$} then due to the "Bochner maximum principle" comes into force [16, p. 30] according to which a subharmonic function on a compact manifold is a constant. As a result, from our Lemma 2 we obtain Corollary 3 which is essentially the Berger-Ebin theorem [2, p. 436]; [29].

\noindent Let us prove our Theorem 3. Let \textit{$\left({\rm M,g}\right)$ }be a complete noncompact Riemannian manifold with nonnegative sectional curvature and ${\rm T}$ be a globally defined non-zero Codazzi 2-tensor with a constant trace. In this case ${\rm Q}_{2} \left(T,\, T\right)\ge 0$ and according to the Theorem 1 the norm $\, \left\| \, {\rm T}\, \right\| $ is a subharmonic function. Therefore, due to the Theorem 1 from [23] we come to the conclusion that there is no non-zero Codazzi tensor ${\rm T}$ on \textit{$\left({\rm M,g}\right)$} such that $\int _{{\rm M}}\left\| \, T\, \right\| \, dVol_{g} <+\infty  $. 

\noindent Let us turn to the Yau theorem [22, p. 663]. Yau's theorem states the following: If $\int _{{\rm M}}\left\| \, {\rm T}\, \right\| ^{{\rm q}} dVol_{g} <+\infty  $\textit{ }for some $q>1$ on a complete \textit{$\left(M,g\right)$ }then $\left\| \, {\rm T}\, \right\| =C$ ($C$ is a constant). In this case, the condition $\int _{{\rm M}}\left\| \, {\rm T}\, \right\| ^{{\rm q}} dVol_{g} <+\infty  $ is equivalent to $C^{q} \int _{{\rm M}}dVol_{g} <+\infty  $. The latter inequality is not feasible on a complete noncompact Riemannian manifold \textit{$\left({\rm M,g}\right)$ }with nonnegative sectional curvature [23]. So $\, \left\| \, T\, \right\| ={\rm C}=0$ that contradicts the existence of a non-zero Codazzi tensor. The Theorem 2 is proved.

\noindent The Corollary 5 does not require any proof. 

\noindent  And now we prove Theorem 4. Let, as before, \textit{$\left({\rm M,g}\right)$} be an ${\rm n}$-dimensional $(n\ge 3)$ locally conformally flat Riemannian manifold with positive constant scalar curvature $s={\rm trace}_{{\rm g}} {\rm Ric}$. In this case, for the traceless Ricci tensor $\overline{Ric}=Ric-n^{-1} s\, g$ we have 
\[Q_{2} \left(\, Ric,\, Ric\right)=Q_{2} \left(\, \overline{Ric,}\, \overline{Ric}\, \right)=\] 
\[=R_{ij} \bar{R}^{ik} \bar{R}_{k}^{j} -R_{ijkl} \bar{R}^{ik} \bar{R}^{jl} =\frac{1}{n-1} \, s\, \left\| \, \overline{Ric}\, \right\| ^{2} +\frac{n}{n-2} \bar{R}_{ij} \bar{R}^{ik} \bar{R}_{k}^{j} \ge \] 
\begin{equation} \label{GrindEQ__3_7_} 
\ge \frac{1}{n-1} \, \left\| \, \overline{Ric}\, \right\| ^{2} \left(s-\sqrt{n\left(n-1\right)} \, \left\| \, \overline{Ric}\, \right\| \, \right) 
\end{equation} 
where, due to Lemma 2.1 from  [47], we used the inequality
\[\bar{R}_{ij} \bar{R}^{ik} \bar{R}_{k}^{j} \ge -\frac{n-2}{\sqrt{n\left(n-1\right)} } \, \left\| \, \overline{Ric}\, \right\| ^{3} .\] 
Then from the formula \eqref{GrindEQ__3_4_} we obtain the inequality 
\begin{equation} \label{GrindEQ__3_8_} 
\frac{1}{2} \, \Delta _{{\rm B}} \, \left\| \, \overline{Ric}\, \right\| ^{2} \ge \frac{1}{n-1} \, \left\| \, \overline{Ric}\, \right\| ^{2} \left(s-\sqrt{n\left(n-1\right)} \, \left\| \, \overline{Ric}\, \right\| \, \right)+\left\| \, \nabla \, \overline{Ric}\, \right\| ^{2} .                 
\end{equation} 
One can prove that $\Delta _{B} \, \left\| \, Ric\, \right\| ^{2} =\Delta _{{\rm B}} \, \left\| \, \overline{Ric}\, \right\| ^{2} $. Therefore, we can rewrite \eqref{GrindEQ__3_8_} in the following form
\begin{equation} \label{GrindEQ__3_9_} 
\frac{1}{2} \, \Delta _{B} \, \left\| \, {\rm Ric}\, \right\| ^{2} \ge \frac{1}{n-1} \, \left\| \, \overline{Ric}\, \right\| ^{2} \left(s-\sqrt{n\left(n-1\right)} \, \left\| \, \overline{Ric}\, \right\| \, \right) 
\end{equation} 
If we suppose that $\, \left\| \, {\rm Ric}\, \right\| ^{2} <\left(n-1\right)^{-1} s^{2} $ then from \eqref{GrindEQ__3_9_} we conclude that $\Delta _{B} \, \left\| \, Ric\, \right\| ^{2} \ge 0$, i.e.$\, \left\| \, Ric\, \right\| ^{2} $ is a nonnegative subharmonic function. 

\noindent In the first case, if \textit{$\left({\rm M,g}\right)$} is a connected complete manifold and\textit{$\, \left\| \, {\rm Ric}\, \right\| ^{2} $} has a global maximum point, then due the ``Hoph maximum principle'' we obtain $\, \left\| \, Ric\, \right\| ^{2} ={\rm C}$ where \textit{C} is a constant. 

\noindent In the second case, we rewrite \eqref{GrindEQ__3_8_} in the following form
\[\left\| \, {\rm Ric}\, \right\| \, \Delta _{B} \, \left\| \, {\it Ric}\, \right\| =\frac{1}{n-1} \, \left\| \, \overline{Ric}\, \right\| ^{2} \left(s-\sqrt{n\left(n-1\right)} \, \left\| \, \overline{Ric}\, \right\| \, \right)+\left\| \, \nabla \, Ric\, \right\| ^{2} -\left\| \, {\rm d}\, \left\| \, {\it Ric}\, \right\| \, \right\| ^{2} \ge \] 
\begin{equation} \label{GrindEQ__3_10_} 
\ge \frac{1}{n-1} \, \left\| \, \overline{Ric}\, \right\| ^{2} \left(s-\sqrt{n\left(n-1\right)} \, \left\| \, \overline{Ric}\, \right\| \, \right) 
\end{equation} 
where we used the first Kato inequality $\left\| \, \nabla \, {\rm Ric}\, \right\| ^{2} \ge \left\| \, d\, \left\| \, {\rm Ric}\, \right\| \, \right\| ^{2} $(see [46]). Let $\, \left\| \, {\rm Ric}\, \right\| ^{2} <\left(n-1\right)^{-1} s^{2} $ then from \eqref{GrindEQ__3_10_} we conclude that $\left\| \, {\rm Ric}\, \right\| $ is a nonnegative subharmonic function. Then for an arbitrary ${\rm q}\ge 2$, either $\int _{{\rm M}}\left\| \, {\rm Ric}\, \right\| ^{{\rm q}} \, dVol_{g} =+\infty  $ or $\left\| \, {\rm Ric}\, \right\| =C$ [22, p. 664]. Therefore, if we suppose that $\int _{{\rm M}}\left\| \, {\rm Ric}\, \right\| ^{{\rm q}} \, dVol_{g} <+\infty  $ at least for one $q\ge 2$, then $\left\| \, {\rm Ric}\, \right\| ={\rm C}$.

\noindent In the third case, we remind that a complete manifold \textit{$\left({\rm M,g}\right)$ }is said to be \textit{parabolic} if it does not admit a positive \textit{Green's function}. In addition, if \textit{$\left({\rm M,g}\right)$} is a parabolic manifold then every subharmonic and bounded-above function on $M$ is constant [20, p. 147]. In our case, $\, \left\| \, Ric\, \right\| ^{2} $ is a subharmonic function on \textit{$\left({\rm M,g}\right)$ }such as $\left\| \, {\rm Ric}\, \right\| ^{2} <\left(n-1\right)^{-1} s^{2} $ for the constant scalar curvature$s>0$. Therefore, for the case of parabolic manifold \textit{$\left({\rm M,g}\right)$ }we have $\, \left\| \, Ric\, \right\| ^{2} ={\rm C}$. 

\noindent Finally, if $\, \left\| \, Ric\, \right\| ^{2} ={\rm C}$ and $\left\| \, {\rm Ric}\, \right\| ^{2} <\left(n-1\right)^{-1} s^{2} $ then from \eqref{GrindEQ__3_9_} we obtain $\, \left\| \, \overline{{\rm Ric}}\, \right\| ^{2} =0$, that is why ${\rm g}$ is an Einstein metric. Therefore, \textit{$\left({\rm M,g}\right)$} becomes a complete Riemannian manifold with positive constant curvature, that means \textit{$\left({\rm M,g}\right)$} is a \textit{spherical space form} [48, p. 69].

\noindent Let us prove our Theorem 5. In the first case, if \textit{$\left\| \, {\rm S}\, \right\| ^{2} \le {\rm n}$} then from \eqref{GrindEQ__2_3_} we conclude that $\left\| \, {\rm S}\, \right\| ^{2} $ is a nonnegative subharmonic function. If, moreover, \textit{$\left({\rm M,g}\right)$} is a connected complete manifold and\textit{$\, \left\| \, {\rm S}\, \right\| ^{2} $} has a global maximum point, then due to the ``Hoph maximum principle'' we obtain $\, \left\| \, Ric\, \right\| ^{2} ={\rm C}$ where \textit{C} is a constant. 

\noindent In the second case, we rewrite \eqref{GrindEQ__2_3_} in the following form
\begin{equation} \label{GrindEQ__3_11_} 
\left\| \, {\it S}\, \right\| \, \Delta _{B} \, \left\| \, {\it S}\, \right\| \ge \left\| \, {\rm S}\, \right\| ^{2} \left(\, n-\left\| \, {\rm S}\, \right\| ^{2} \right) 
\end{equation} 
Let \textit{$\left\| \, {\rm S}\, \right\| ^{2} \le {\rm n}$} then from \eqref{GrindEQ__3_11_} we conclude that $\left\| \, {\rm S}\, \right\| $ is a nonnegative subharmonic function. Then for an arbitrary ${\rm q}\ge 1$, either $\int _{{\rm M}}\left\| \, {\rm S}\, \right\| ^{{\rm q}} \, dVol_{g} =+\infty  $ or $\left\| \, {\rm S}\, \right\| =C$ [22, p. 664]. Therefore, if we supouse  $\int _{{\rm M}}\left\| \, {\rm S}\, \right\| ^{{\rm q}} \, dVol_{g} <+\infty  $ for some $q\ge 2$, then $\left\| \, {\rm S}\, \right\| =C$.

\noindent In the third case, if \textit{$\left({\rm M,g}\right)$} is a parabolic manifold then \textit{$\left\| \, {\rm S}\, \right\| ^{2} ={\rm C}$} because \textit{$\, \left\| \, {\rm S}\, \right\| ^{2} $} is a subharmonic function on \textit{$\left({\rm M,g}\right)$ }such that \textit{$\left\| \, {\rm S}\, \right\| ^{2} \le {\rm n}$}. 

\noindent Finally, if  \textit{$\left\| \, {\rm S}\, \right\| ^{2} ={\rm C}$} and \textit{$\left\| \, {\rm S}\, \right\| ^{2} \le {\rm n}$ }then from \eqref{GrindEQ__2_3_} we obtain either $\left\| \, {\rm S}\, \right\| ^{2} \equiv 0$ or $\left\| \, {\rm S}\, \right\| ^{2} \equiv n$.\textit{ } In the first case, $\left({\rm M,}\, {\rm g}\right)$ must be totally geodesic. At the same time, we know that an arbitrary ${\rm n}$-dimensional complete totally geodesic submanifold of the sphere ${\rm S}^{n+p} $ is a sphere ${\rm S}^{n} $ [49]. Therefore, $\left({\rm M,}\, {\rm g}\right)$ has to be an equator ${\rm S}^{n} \subset {\rm S}^{n+1} $. In the second case, $\left({\rm M,}\, {\rm g}\right)$ is locally isometric to a generalized Clifford torus ${\rm S}^{{\rm k}} \left(r_{1} \right)\times S^{n-k} \left(r_{2} \right)$, which is the standard product embedding of the product of two spheres of radius ${\rm r}_{1} =\sqrt{k\, n^{-1} } $ and $r_{2} =\sqrt{\left(n-k\right)\, n^{-1} } $, respectively [50].

\bigskip
\noindent {\footnotesize \textbf{Acknowledgments. }Our work was supported by the Russian Foundation for Basic Research of the  Russian  Academy of Science (projects Nos. 16-01-00053 and 16-01-00756). The authors would like to express their thanks to Professor Nelli Makievskaya for her editorial suggestions, that contributed to the improvement of this paper as well as to the referee for his careful reading of the previous version of the manuscript and suggestions about this paper which led to various improvements.}

\noindent \textbf{References}

\noindent [1] Derdzinski A., Shen C., ``Codazzi tensor fields, curvature and Pontryagin forms``, Proc.  London Math. Soc., \textbf{47}:3 (1983), 15-26.

\noindent [2] Besse, A. Einstein Manifolds, Springer-Verlag, Berlin-Heidelberg (1987).

\noindent [3] Simon U., Schwenk-Schellschmidt A., Viesel H., Introduction to the affine differential geometry of hypersurfaces, Science University Tokyo Press, Tokyo (1991).

\noindent [4] Catino G., Mantegazza C., Mazzieri L., ``A note on Codazzi tensors'', \textit{Mathematische Annalen}, \textbf{362}: 1-2 (2015), 629-638.

\noindent [5] Merton G., ``Codazzi tensors with two eigenvalue functions'', \textit{Proc. Amer. Math. Soc.}, \textbf{141} (2013), 3265-3273.

\noindent [6] Aledo J.A., Espinar J.M., Galvez J.A., ``The Codazzi equation for surfaces'', \textit{Advances in Math}., \textbf{224} (2010), 2511-2530.

\noindent [7] Bonsante F., Seppi A., ``On Codazzi Tensors on a hyperbolic surface and flat Lorentzian geometry'', \textit{International Mathematics Research Notesá} \textbf{2016}: 2 (2016), 343-417.

\noindent [8] Catino G., Mantegazza C., Mazzieri L., ``A note on Codazzi tensors'', \textit{Mathematische Annalen}, \textbf{362}:1-2 (2015), 629-638.

\noindent [9] Baum H., ``Holonomy group of Lorentzian manifolds: A status report'', \textit{Global differential geometry},  Springer-Verlag,  Berlin-Heidelberg (2012), 163-200.

\noindent [10] Catino C., ``On conformally flat manifolds with constant positive scalar curvature'', \textit{Proc. Amer. Math. Soc}., \textbf{144} (2016), 2627-2634.

\noindent [11] Stepanov S. E., "Fields of symmetric tensors on a compact Riemannian manifold'', \textit{Math. Notes}, \textbf{52}:4 (1992), 1048-1050.

\noindent [12] Liu H.L., Simon U., Wang C.P., ``Codazzi tensor and the topology of surfaces'', \textit{Annals of Global Analysis and Geometry}, \textbf{16} (1998), 189-202.

\noindent [13] Liu H.L., Simon U., Wang C.P., ``Higher order Codazzi tensors on conformally flat spaces'', \textit{Beitr\"{a}ge zur Algebra und Geometrie~Contributions to Algebra and Geometry},~ \textbf{39}:2 (1998) 329-348.

\noindent [14] Leder J., Schwenk-Schellschmidt A., Simon U., Wiehe M., ``Generating higher order Codazzi tensors by functions''\textit{, Geometry and Topology of Submanifolds }IX'', London and Singapore, Word Scientific Publishing Co. Pte. Ltd. (1999),  174-191.

\noindent  [15] Huang G., ``Rigidity of Riemannian manifolds with positive scalar curvature'', Annals of Global Analysis and Geometry, 54:2 (2018), 257-272.

\noindent [16] Yano K., Bochner S., Curvature and Betti numbers,  Princeton, Princeton Univ. Press (1953).

\noindent [17] Wu H., The Bochner technique in differential geometry,  New York, Harwood Academic Publishers (1988).

\noindent [18] Peterson P., Riemannian geometry, Springer Int. Publ. AG, Switzerland   (2016).

\noindent [19] B\'{e}rard P.H., ``From vanishing theorems to estimating theorems: the Bochner technique revisited'', \textit{Bulletin of the American Math. Soc}., \textbf{19}:2 (1988), 371-406.

\noindent [20] Pigola S., Rigoli M., Setti A.G., Vanishing and Finiteness Results in Geometric Analysis. A Generalization of the Bochner Technique, Birkh\"{a}user Verlag AG,  Berlin (2008).

\noindent [21] Calabi E., ``An extension of E. Hopf's maximum principle with an application    to Riemannian geometry'', \textit{Duke Math. J.}, \textbf{25} (1958), 45-56.

\noindent [22] Yau S.T., ``Some function-theoretic properties of complete Riemannian manifold and their applications to geometry'', \textit{ Indiana  Univ. Math. J}., \textbf{25}:7\textbf{ }(1976), 659-679.

\noindent [23] Greene R.E., Wu H., ``Integrals of subharmonic functions on manifolds of nonnegative curvatures'', \textit{Inventions Math}., \textbf{27} (1974), 265-298.

\noindent [24]  Li P., Schoen R., ``${\rm L}^{{\rm p}} $and mean value properties of subharmonic functions on Riemannian manifolds'', \textit{Acta Mathematica}, \textbf{153 }(1984), 279-301.

\noindent [25] Stepanov S.E., Tsyganok I.I., ``Conformal Killing \textit{L}${}^{2}$-forms on complete Riemannian manifolds with nonpositive curvature operator'', \textit{Journal of Mathematical Analysis and Applications}, \textbf{458}:1 (2018), 1-8.

\noindent [26] Stepanov S.E., Tsyganok I.I., ``Theorems on conformal mappings of complete Riemannian manifolds and their applications'', \textit{Balkan Journal of Geometry and its Applications}, \textbf{22}:1 (2017), 81-86. 

\noindent [27] Stepanov S.E., Mike\v{s} J., ``Liouville-type theorems for some classes of Riemannian almost product manifolds and for special mappings of Riemannian manifolds'', Differential Geometry and its Applications, 54 (2017), Part A, 111-121.

\noindent [28] Tani M., ``On a conformally flat Riemannian space with positive Ricci curvature'', \textit{Tohoku Math. Journal}, \textbf{19}:2 (1967), 227-231.

\noindent [29] Berger M., Ebin D., ``Some decomposition of the space of symmetric tensors on a Riemannian manifold'', \textit{Journal of Differential Geometry}, \textbf{3 }(1969), 379-392.

\noindent [30] Goldberg S.I., ``An application of Yau's maximum principle to conformally flat spaces'', \textit{Proceeding of the American Mathematical Society}, \textbf{79}:2 (1980), 260-270.

\noindent [31] Bouguignon J.-P., H. Karcher, ``Curvature operators: pinching estimates and geometric examples'', \textit{Ann. Sc. \'{E}c. Norm. Sup.}, \textbf{11} (1978), 71-92.

\noindent [32] Tachibana S., Ogiue K., ``Les vari\'{e}t\'{e}s riemanniennes dont l'op\'{e}rateur de coubure restreint est positif sont des sph\`{e}res d'homologie r\'{e}elle'', \textit{C. R. Acad. Sc.  Paris}, \textbf{289} (1979), 29-30.

\noindent [33]  Stepanov S. E., Mikes. J, "Hodge-de Rham Laplacian and Tachibana operator on a compact Riemannian manifold with a curvature operator of fixed sign'', \textit{Izvestiya}: \textit{Mathematics}, \textbf{79}:2 (2015) 375-387. 

\noindent [34] Kashiwada T., ``On the curvature operator of the second kind'', \textit{ Natural  Science  Report} \textit{ Ochanomizu  University}, \textbf{44}:2 (1993), 69-73.

\noindent [35] Stepanov S.E., Tsyganok I.I., ``Theorems of existence and of vanishing of conformally Killing forms'', \textit{Russian Mathematics, }\textbf{58}:10 (2014), 46-51.

\noindent [36] Stepanov S.,  Tsyganok I.I., ``Vanishing theorems for harmonic mappings into non-negatively curved manifolds and their applications'', \textit{Manuscripta Mathematica}, 154:1-2 (2017), 79-90. 

\noindent [37] Bettiol R.G., Mendes R.A.E., ``Sectional curvature and Weitzenb\"{o}ck formulae'', \textit{Preprint}, arXiv:1708.09033v1 [math.DG] 29 Aug 2017, 24 pp.

\noindent [38] Guan P., Viaclovsky J., Wang G., ``Some properties of the Schouten tensor and applications to conformal geometry'', \textit{Transactions of the American Mathematical Society}, \textbf{355} (2002), no. 3, 925-933.

\noindent [39] Bourguignon J.P., ``Formules de Weitzenb\"{o}k en dimension 4'', Seminare A. Besse sur la geometrie Riemannienne dimension 4, Cedic. Ferman, Paris, (1981), 156-177.

\noindent [40] Bourguignon J.P., ``Les vari\'{e}t\'{e}s de dimension 4 \'{a} signature non nulle dont la courbure est harmonique sont d'Einstein'', \textit{Invent. Math., }63 (1981), 263-286.

\noindent [41] Eisenhart L. P., Riemannian geometry, Princeton Univ. Press (1926).

\noindent [42] Kobayashi Sh., Nomizu K., Foundations of differential geometry, vol. 2, Interscience Publishers, New York-London-Sydney (1969).

\noindent [43] Donaldson S.K., ``Symmetric spaces, K\"{a}hler geometry and Hamiltonian dynamics'', \textit{Amer. Math. Soc. Transl}. \textbf{196} (1999), 13-33.

\noindent [44] Wegner B., ``Kennzeichnungen von R\"{a}umen konstanter Kr\"{u}mmung unter lokal konform Euklidischen Riemannschen Mannigfaltigkeiten``, \textit{Geom. Dedicata}, \textbf{2} (1973) 269--281.

\noindent [45] Peng C.-K., Terng C.-L., ``The scalar curvature of minimal hypersurfaces in spheres'', \textit{Math. Ann}., \textbf{266} (1983), 105-113.

\noindent [46] Calderbank D.M.J., Gauduchon P., Herzlich M., ``Refined Kato inequalities and conformal weights in Riemannian geometry'', \textit{Journal of Functional Analysis}, \textbf{173} (2000), 214-255.

\noindent [47] Okumura M., ``Hypersurfaces and a pinching problem on the second fundamental tensor'', \textit{Amer. J. Math.}, \textbf{96} (1974), 207--213. 

\noindent [48] Wolf G., Spaces of constant curvature,  California Univ. Press,  Berkley (1972).

\noindent [49] Morvan J.-M., ``Differential geometry of Riemannian submanifolds: recent results'', Proceedings of the Conference on Algebra and Geometry ( Kuwait, 1981), Kuwait Univ. Press,  Kuwait (1981), 145-151.

\noindent [50] Lawson B., ``Local rigidity theorems for minimal hypersurfaces'', \textit{Ann. of Math}., \textbf{89}:2 (1969), 187-197.

\noindent

\end{document}